\documentclass[reqno]{amsart}
\usepackage[greek,english]{babel}
\usepackage[iso-8859-7]{inputenc}
\usepackage{amsthm}
\usepackage{amsfonts}
\usepackage{graphicx}
\newtheorem{thm}{Theorem}

\newtheorem{lem}{Lemma}
\newtheorem{prop}{Proposition}
\newtheorem{rem}{Remark}
\begin{document}

\title[Kinks in an array of coupled pendula]
 {Traveling kinks in an infinite array of weakly coupled pendula}

\author{ Christos Sourdis}

\address{University of Athens, Greece.}

\email{sourdis@uoc.gr}

\begin{abstract}
We prove the existence of heteroclinic traveling waves (kinks) in an infinite array of weakly coupled pendula. Our approach is to apply a perturbation argument from
the anti-continuum  limit.
\end{abstract}

\maketitle

\section*{Introduction and main result}

\subsection*{The model}

The following infinite system of linearly coupled ODEs describes the motion of an
array of pendula each of which is coupled to its nearest neighbors
by a torsional spring with a coupling coefficient $k$:

\begin{equation}\label{eq----1}
\ddot{v}_n- k (v_{n+1}-2v_n+v_{n-1})+\sin(v_n)=0,\ \
\ n \in \mathbb{Z}.
\end{equation}
More precisely,  $v_i$ represents the angle  formed by the $i$th pendulum
with the vertical axis (assuming physical
units have been scaled appropriately).
We refer the reader to \cite{kresetrudi,levi,scott} for more
details on the physical background of the problem.

The system (\ref{eq----1}) is known as the \emph{discrete sine-Gordon equation} and also serves as a
model of  arrays of Josephson
junctions \cite{imry}, or as a dynamical Frenkel-Kontorova model of
electrons in a crystal lattice \cite{braunK}.

\subsection*{Traveling wave solutions}
We shall construct solutions of the above equation in the form of
traveling waves (cf. \cite{chowMallet}). To this end, we let $z=n$, write $v_n(t)=v(z,t)$, and seek a
solution of the form $v(z,t)=v(\xi)=v(z-ct)$ satisfying the equation
\begin{equation}\label{eq0}
 c^2
\frac{d^2v}{dz^2}-k[v(z+1)-2v(z)+v(z-1)]+\sin(v)=0.
\end{equation}
\subsubsection*{Known results}
In \cite{Schmidt} explicit kink solutions to the above equation were given for a certain nonlinearity $f$ (in place of the sine) which satisfies (H1) below.   Periodic traveling wave solutions to (\ref{eq----1}) have been shown to exist recently in \cite{rothos,rothos2,levia} in the strong coupling regime (i.e. when $k\gg 1$). This was achieved in \cite{rothos} using techniques from dynamical systems; in \cite{rothos2} using variational and topological techniques; in \cite{levia} using a fixed point argument.
However, as is pointed out in \cite{rothos}, kink solutions to (\ref{eq0}) connecting $-\pi$ to $\pi$ should not exist in this regime.
Nevertheless, such solutions were constructed variationally in \cite{kreiner} provided that $c$ is sufficiently large. The question of persistence of kink solutions in the continuum limit of (\ref{eq----1}) and (\ref{eq4model}) was discussed in \cite{aigner,dmitriev,IossPeli,oxt,savin}.
We also refer to \cite{rothos3,kouk} for further results on the existence of localized structures in long range interaction lattices (stationary or traveling).
\subsection*{The main result}\label{sec2}
In the current paper, for given $c \neq 0$, we study heteroclinic
waves for sufficiently small $k$, which from  now on we will call $\varepsilon$.
We will also consider a more general class of nonlinearities
 $f$ satisfying (H1) below, covering both the discrete sine-Gordon equation as well as the important $\phi^4$ model \cite{BeK}:
 \begin{equation}\label{eq4model}
\ddot{v}_n- k (v_{n+1}-2v_n+v_{n-1})+2(v_n-v_n^3)=0,\ \
\ n \in \mathbb{Z},
\end{equation}
see Remark \ref{rem} below.
 Moreover, motivated from \cite{bateschenchmaj}, we will allow infinite range
and not just nearest neighbor or finite length interaction, although
those are included as special cases.

The equation we will be dealing with  is
\begin{equation}\label{eq1}
u''-\epsilon \sum_{k=-\infty}^{\infty}a_ku(z-k)+f(u)=0,\ \ \ z \in
\mathbb{R},
\end{equation}
together with the conditions
\begin{equation}\label{eq2}
\lim_{z \to -\infty}u(z)=0,\ \ \ \lim_{z \to \infty}u(z)=1.
\end{equation}
Here $\epsilon \geq 0$ is small, and we assume that

\textbf{(H1)} $f \in C^2(\mathbb{R}),\ f(0)=f(1)=0,\ \
f'(0),\ f'(1)<0;$
\[
F(1)=0,\ F(u)\neq 0\ \ \forall u \in (0,1)
\]
where $F(u)=\int_{0}^{u}f(s)ds.$
\\

\textbf{(H2)} $\sum_{k=-\infty}^{\infty}a_k=0,\ a_0<0,\ a_k=a_{-k},\
\textrm{and}\ \sum_{k \geq 1}^{}|a_k|k^2 <\infty.$
\\
When $\epsilon=0$, that is in the so called anti-continuum limit \cite{makay}, equation (\ref{eq1}) becomes
\begin{equation}\label{eq3}
u''+f(u)=0.\end{equation} Under the hypotheses (H1), the above
equation has a heteroclinic solution $u_0$ satisfying (\ref{eq2})
(see \cite{arnold}).

Our result is
\begin{thm}\label{thm1}
If $\epsilon>0$ is sufficiently small, then there exists a solution
$u_\epsilon$ of (\ref{eq1}) such that
\[
||u_\epsilon-u_0||_{H^2(\mathbb{R})}\leq C\epsilon
\]
($C>0$ is a constant independent of $\epsilon$).
\end{thm}
\begin{rem}\label{rem}
The choice of the roots of $f$ to be $0$ and $1$ is made for convenience purposes only and causes no loss of generality. For instance, the traveling kink problem (\ref{eq0}) can be embedded in our framework  by plainly letting
\[
u=\frac{v+\pi}{2\pi}\ \ \textrm{and}\ \ f(u)=\frac{1}{2\pi}\sin(2\pi u-\pi).
\]
Similarly, the corresponding change of variables for (\ref{eq4model}) is
\[
u=\frac{v+1}{2}\ \ \textrm{and}\ \ f(u)=(2u-1)-(2u-1)^3 ,
\]
(keep in mind the first assumption in (H2)).
\end{rem}

\subsubsection*{Method of proof}
To prove this we adapt a
 technique from an earlier paper of ours (cf. \cite{AFFS}, but see also \cite{dPK} for a related idea). We use
 two important properties:

 (i)\ \underline{Nondegeneracy of $u_0$} The operator obtained by linearizing the
left-hand side of (\ref{eq3}) at $u_0$ has $0$ as a simple isolated
eigenvalue, the remaining spectrum being in the open left half-plane
(\cite{henry}, Section 5.4).

(ii)\ \underline{Hamiltonian form of the problem} The equation (\ref{eq1}) arises when seeking traveling waves to the lattice equation
\[
 \ddot{u}_n-\epsilon\sum_{k=-\infty}^{\infty}a_ku_{n-k}+f(u_n)=0,\ \ n\in \mathbb{Z},
 \]
which comes from the Hamiltonian on $\ell^2\times \ell^2$ defined by
\[
H(\textbf{p},\textbf{u})=\sum_{n}^{}\left(\frac{1}{2}p_n^2+\epsilon \sum_{m}^{}a_{n-m}(u_n-u_m)^2+F(u_n) \right),
\]
where $\textbf{p}=(\dot{u}_n)\in \ell^2$ and $\textbf{u}=({u}_n)\in \ell^2$.

As is explained in \cite{shatah}, the persistence of heteroclinic or homoclinic orbits of Hamiltonian systems under Hamiltonian perturbations is
a delicate issue, as the Melnikov integral vanishes. An analogous difficulty also arises in the lattice setting at hand, as we will discuss  in Section \ref{seclyap}.
\subsection*{Outline of the paper}
 In Section \ref{sec3} we present the proof of Theorem \ref{thm1}, in
 Section \ref{seclyap} we make some comments on the application of the standard Lyapunov-Schmidt reduction,
 and in the Appendix's we prove some technical lemmas.
\subsection*{Notation}
In what follows, $||\cdot||_{L^2}$, $||\cdot||_{L^\infty}$, and
$||\cdot||_{H^i},\ (i=1,2)$ denote the norms of the spaces
$L^2(\mathbb{R})$, $L^\infty(\mathbb{R})$, and $H^i(\mathbb{R})$,
respectively. Also,
\[
(\phi,\psi)\equiv \int_{\mathbb{R}}^{}\phi \psi dz;\ \ \ \
\phi\bot\psi  \  \Leftrightarrow \  (\phi,\psi)=0.
\]
Unless specified otherwise $C/c$ denotes a large/small positive
constant independent of $\epsilon>0$ whose value will change from
line to line. In many cases we will not explicitly write the obvious
dependence of functions on $\epsilon$.

\section{Proof of Theorem \ref{thm1}}\label{sec3}
\subsection{Properties of the linear operator $\Delta$}\label{sec4}
Consider the linear operator $\Delta$ defined via
\[
\Delta u := \sum_{k=-\infty}^{\infty}a_ku(z-k).
\]
Then one can verify that if $u,v \in L^2(\mathbb{R})\cap
L^\infty(\mathbb{R})\cap C(\mathbb{R})$, then
\begin{equation}\label{eq4}
||\Delta u||_{L^2}\leq C ||u||_{L^2}\ \ \ (C\ \textrm{independent of
}\ u),
\end{equation}
and
\[
(\Delta u,v)=(u,\Delta v).
\]
We also have the following lemma whose proof will be given in
Appendix \ref{apen1}.
\begin{lem}\label{lem1}
If $u \in C^1(\mathbb{R})$, $u'\in L^2(\mathbb{R})$;
\begin{equation}\label{eqlimits}
\lim_{z\to -\infty}u(z)=u(-\infty)\in \mathbb{R} \ \ \ \textrm{and}
\ \ \ \lim_{z \to \infty}u(z)=u(\infty)\in \mathbb{R},
\end{equation}
\textrm{then}
\[ \Delta u \in L^2(\mathbb{R}) \ \ \ and \ \ \
(\Delta u,u')=0.\]
\end{lem}
\subsection{Properties of the heteroclinic $u_0$}
It is well known (see \cite{arnold}, \cite{halekocak}) that $u_0$ is
the unique (up to translation) solution of (\ref{eq3}), (\ref{eq2}).
Furthermore $u_0'>0$, $u_0$ approaches its limits exponentially and
\[
u_0'(z),\ |u_0''(z)|,\ |u_0'''(z)|\leq Ce^{-c|z|},\ \ \ z \in
\mathbb{R}.
\]
The linear operator $L^0$ with $D(L^0)=H^2(\mathbb{R})$ and
\[
L^0\phi=\phi'' +f'(u_0(z))\phi
\]
is self-adjoint in $L^2(\mathbb{R})$ and
$\sigma(-L^0)\subseteq\{0\}\cup [c,\infty)$ with $0$ a simple
eigenvalue corresponding to $u_0'$. (Note that since $f'(0),\
f'(1)<0$ then $\sigma_{ess}(-L^0)\subseteq [c,\infty)$; thus the
only thing left to prove is that $0$ is the principal eigenvalue
which follows from $u_0'>0$.)

These properties imply the following important proposition.
\begin{prop}\label{prop1}
Let $g \in L^2(\mathbb{R})$ with $g\bot u_0'$, then there exists a
unique $\phi \in H^2(\mathbb{R})$ with $\phi\bot u_0'$ such that
\[
L^0\phi=g.
\]
Furthermore, we have
\[
||\phi||_{H^2}\leq C ||g||_{L^2}
\]
with $C$ independent of g.
\end{prop}
\subsection{The perturbation argument} We search for a solution of
(\ref{eq1}) in the form $u_\epsilon=u_0+\phi_\epsilon$ with
$\phi_\epsilon \in H^2(\mathbb{R})$ and $\phi_\epsilon \bot u_0'.$
Then, the fluctuation $\phi_\epsilon$ must satisfy
\[
L^0\phi=\epsilon \Delta \phi +N(\phi)+\epsilon\Delta u_0=E(\phi)
\]
with \[ N(\phi)=-f(u_0+\phi)+f(u_0)+f_u(u_0)\phi.
\]

We note that $\left(E(\phi_\epsilon),u_0' \right)=0$ holds. However, if iterations of the form $L^0\phi_{n+1}=E(\phi_n)$ were to be performed, for capturing the desired $\phi_\epsilon$ in the limit $n\to \infty$, the iteration $\phi_n$ may not satisfy this orthogonality condition which is necessary for solving for $\phi_{n+1}\in H^2(\mathbb{R})$. To deal with this issue, at each step of the iteration we will project $E(\phi_n)$ to $\{u_0'\}^\bot$ and then solve for the corresponding $\phi_{n+1}$.

We  thus define a mapping ${T}:H^2(\mathbb{R})\cap \{u_0'\}^\bot
\to H^2(\mathbb{R})\cap \{u_0'\}^\bot$ via ${T}(\phi) =
{\psi}$ where
\begin{equation}\label{eq5}
L^0{\psi}=-b(\phi)u_0'+\epsilon \Delta \phi
+N(\phi)+\epsilon\Delta u_0
\end{equation}
and
\begin{equation}\label{eq6}
b(\phi)=\frac{1}{||u_0'||^2_{L^2}}\left(\epsilon \Delta \phi
+N(\phi)+\epsilon \Delta u_0,u_0'\right).
\end{equation}
Note that ${T}$ is well defined via Proposition $\ref{prop1}$
since the right hand side of (\ref{eq5}) is orthogonal to $u_0'$.
(Note also that by Lemma $\ref{lem1}$ we have $\Delta u_0 \in
L^2(\mathbb{R}).$)

Let
\[
B_\epsilon =\{ \phi \in H^2(\mathbb{R})\cap \{ u_0'\} ^ \bot  : \
||\phi ||_{H^2} \leq M \epsilon \}
\]
with $M$ a positive constant independent of $\epsilon>0$ to be
determined later.
We will show that there exists a large $M>0$ such that, provided
$\epsilon >0$ is sufficiently small, ${T}$ maps $B_\epsilon$
into itself and is a contraction.
 Let $\phi \in B_\epsilon$, then via
(\ref{eq5}), (\ref{eq6}) and Proposition \ref{prop1},
\begin{eqnarray}
 ||{\psi}||_{H^2} &\leq& C \left(|b(\phi)|+\epsilon ||\Delta \phi||_{L^2}
+||N(\phi)||_{L^2}+\epsilon ||\Delta u_0||_{L^2}\right)\nonumber
 \\
   & & \label{eq8} \\
   &\stackrel{(\ref{eq6})}{\leq}& C \left(\epsilon ||\Delta \phi||_{L^2}
+||N(\phi)||_{L^2}+\epsilon ||\Delta u_0||_{L^2}\right). \nonumber
\end{eqnarray}
The first and third term will be estimated from (\ref{eq4}) and
Lemma \ref{lem1} respectively. To estimate the nonlinear term
$N(\phi)$, we first recall the embedding $||\phi||_{L^\infty}\leq C
||\phi||_{H^1}$ for every $\phi \in H^1(\mathbb{R})$. Hence, setting
$C=\sup_{|s|\leq 2}|f''(s)|$, we have
\begin{equation}\label{eqN}
|N(\phi)|\leq CM\epsilon|\phi| \ \ \textrm{and} \ \
|N(\phi_1)-N(\phi_2)|\leq CM\epsilon |\phi_1-\phi_2|
\end{equation}
pointwise for all $\phi,\phi_1,\phi_2 \in B_\epsilon$. Thus,
(\ref{eq8}) yields
\begin{eqnarray}
 ||{\psi}||_{H^2} &\leq& C \left(\epsilon || \phi||_{L^2}
+M\epsilon||\phi||_{L^2}+\epsilon \right)\nonumber
 \\
   & & \nonumber \\
   &\leq& C \left(M\epsilon
+M^2\epsilon +1\right)\epsilon. \nonumber
\end{eqnarray}
Choosing a large $M$ (say $M=2C$), then $||{\psi}||_{H^2}\leq
M\epsilon$ provided $\epsilon>0$ is sufficiently small, i.e.
${T}:B_\epsilon \to B_\epsilon$. Similarly we can show that
${T}$ is a contraction in the $H^2$ norm.

Since
$B_\epsilon$ is closed with respect to this norm, the Banach fixed
point theorem gives us a fixed point $\phi_* \in B_\epsilon$ of
${T}$. Then
\begin{equation}\label{eq9}
u_\epsilon=u_0+\phi_*^\epsilon\end{equation} satisfies \begin{equation}\label{eq7}
u''-\epsilon \Delta u +f(u)=-b_\epsilon u_0'\end{equation}
for some $b_\epsilon \in \mathbb{R}$
($b_\epsilon=b(\phi_*^\epsilon)$). Multiplying (\ref{eq7}) by
$u'=u_\epsilon'$ and integrating over $\mathbb{R}$ yields
\[
\int_{\mathbb{R}}^{}u'u''dz-\epsilon(\Delta
u,u')+\int_{\mathbb{R}}^{}f(u)u'dz=-b_\epsilon
\int_{\mathbb{R}}^{}u_0'u'dz.
\]
Since $\phi_*\in H^2(\mathbb{R})$, we have that $u(-\infty)=0$,
$u(\infty)=1$ and $u'(\pm \infty)=0$. The left hand side of the
above equation is 0 (see Lemma \ref{lem1} and recall that
$F(1)=0$) and we get that
\[
b_\epsilon\int_{\mathbb{R}}^{}u_0'u'dz=0.
\]
This implies that $b_\epsilon=0$ since
\begin{eqnarray}
  \int_{\mathbb{R}}^{}u_0'u'dz &=& \int_{\mathbb{R}}^{}u_0'(u_0'+\phi_*')dz \geq
\int_{\mathbb{R}}^{}u_0'^2dz-||u_0'||_{L^2}||\phi_*'||_{L^2}\geq \nonumber\\
    & &  \nonumber\\
   &\geq&\int_{\mathbb{R}}^{}u_0'^2dz- C\epsilon>0,\nonumber
\end{eqnarray}
provided $\epsilon>0$ is sufficiently small.

Therefore $u$ given by (\ref{eq9}) is a solution of (\ref{eq1})
satisfying the estimate of Theorem \ref{thm1}, thereby completing
the proof.

\begin{rem}\label{remODDH}
  If $f(u+1/2)$ is odd, as it would be the case for applications to (\ref{eq0}) and (\ref{eq4model}), then the proof of Theorem \ref{thm1} becomes much simpler. Indeed, we can apply our fixed point argument restricted to the space of functions such that $\phi-1/2$ is odd, without the need to introduce the projection operator. Keep in mind that $u_0-\frac{1}{2}$ would be odd and $u_0'$ would be even. So, the orthogonality condition $\left(E(\phi),u_0' \right)=0$ would be automatically satisfied.
\end{rem}

\begin{rem}
If the $\epsilon=0$ equation (\ref{eq1}) has a unique even
homoclinic solution $u_0$ (see \cite{bl} for necessary and
sufficient conditions on $f$), then the proof of persistence for
$\epsilon$ small is considerably simplified by seeking
$u_\epsilon=u_0+\phi$ with $\phi \in H^2(\mathbb{R})$ even. Note
that given an even $g\in L^2(\mathbb{R})$, there exists a unique
even $\phi \in H^2(\mathbb{R})$ such that $\phi''+f'(u_0)\phi=g$.
Moreover $||\phi||_{H^2}\leq C||g||_{L^2}$ for some $C>0$
independent of $g$. This problem was briefly discussed at the end of \cite{BZ}.
\end{rem}
\section{Some remarks on the standard Lyapunov-Schmidt
approach}\label{seclyap} In this section we make some remarks on a
difficulty that arises when trying to prove Theorem \ref{thm1} using
the standard Lyapunov-Schmidt reduction.

We have seen that $u_0$ satisfies (\ref{eq1}) up to an order of
$\epsilon$. We begin by refining this approximation so that
$u_{ap}=u_0+\epsilon u_1$ satisfies (\ref{eq1}) up to an order of
$\epsilon^2$. We choose $u_1 \in H^2(\mathbb{R})$, $u_1\bot u_0'$
such that
\begin{equation}\label{equ1}
u_1''+f'(u_0)u_1=\Delta u_0
\end{equation}
(this is possible via Lemma \ref{lem1} and Proposition \ref{prop1}).
Then, a simple calculation gives
\[
-G(\epsilon):=u_{ap}''-\epsilon \Delta
u_{ap}+f(u_{ap})=-\epsilon^2\Delta u_1-N(\epsilon u_1)
\]
and thus from (\ref{eqN}):
\begin{equation}\label{eqG}
||G(\epsilon)||_{L^2}\leq C \epsilon^2.
\end{equation}

 We seek a solution of (\ref{eq1}) in the form
$u_\epsilon=u_{ap}+\psi_\epsilon$ with $\psi_\epsilon \in
H^2(\mathbb{R})$. Then $\psi_\epsilon$ must satisfy
\begin{equation}\label{eq10}
L^\epsilon\psi=N_{{ap}}(\psi)+G(\epsilon)
\end{equation}
where $L^\epsilon \psi=\psi''+f'(u_{ap})\psi-\epsilon \Delta \psi$
and $N_{{ap}}(\psi)=-f(u_{ap}+\psi)+f(u_{ap})+f'(u_{ap})\psi$.

Since $\Delta: H^2(\mathbb{R})\to L^2(\mathbb{R})$ is a bounded
linear operator and $||u_{ap}-u_0||_{L^\infty}\leq C\epsilon$,
$L^\epsilon$ is a regular $O(\epsilon)$ perturbation of $L^0$. From
the special form of the perturbation, however, the  simple
eigenvalue $0$ of $L^0$ is perturbed to an $O(\epsilon^2)$ simple
eigenvalue of $L^\epsilon$ (this is the  source of the
difficulty). We point out that such small eigenvalues would not have been present if we were in the symmetric setting of Remark \ref{remODDH}. More precisely we have the following Proposition whose
proof we postpone to Appendix \ref{apen2}.
\begin{prop}\label{prop2}
If $\epsilon \geq 0$ is sufficiently small then $\sigma(-L^\epsilon)
\subset \{\lambda_1(\epsilon) \}\cup [c,\infty)$ with
$\lambda_1(\epsilon)$ simple corresponding to $\phi_1(\epsilon)\in
H^2(\mathbb{R})$ with $||\phi_1(\epsilon)||_{H^2}=1$. Moreover
$\lambda_1(\epsilon)$, $\phi_1(\epsilon)$ depend smoothly on
$\epsilon$ up to $\epsilon=0$ and
\begin{equation}\label{Lepsilon}
\begin{array}{ll}
  \lambda_1(\epsilon)= O(\epsilon^2) \\
   &  \\
\phi_1(\epsilon)=    \frac{u_0'}{||u_0'||_{H^2}}+O(\epsilon)\ \ \ \
(\textrm{here}\ ||O(\epsilon)||_{H^2}\leq C \epsilon ).
\end{array}
\end{equation}
\end{prop}
Define the orthogonal projection $P$ onto the span of $\phi_1$ by
\[
P\psi=(\psi,\phi_1(\epsilon))\frac{\phi_1(\epsilon)}{||\phi_1(\epsilon)||^2_{L^2}}.
\]
According to this projection we have
\[
H^2(\mathbb{R})=\textrm{span}\{\phi_1\}\oplus X_1,\ \ \
L^2(\mathbb{R})=\textrm{span}\{\phi_1\}\oplus Y_1,
\]
where $X_1,\ Y_1$ are respectively the kernel of $P$ in
$H^2(\mathbb{R})$ and $L^2(\mathbb{R})$. By decomposing $\psi$ as
$\psi=a\phi_1(\epsilon)+v$ $(a\in \mathbb{R}, v \in X_1)$, one finds
that (\ref{eq10}) is equivalent to
\begin{equation}\label{eq11}
\begin{array}{rl}
  L^\epsilon v= & (I-P)\{N_{ap}(a\phi_1(\epsilon)+v)+G(\epsilon)\}  \\
  &  \\
  -a\lambda_1(\epsilon)\phi_1(\epsilon)= & P\{N_{ap}(a\phi_1(\epsilon)+v)+G(\epsilon)\}.
\end{array}
\end{equation}
Applying Proposition \ref{prop2}, using (\ref{eqG}), and the Banach
fixed point theorem, we can uniquely solve $(\ref{eq11})_{(i)}$ for
$v=v^*(a,\epsilon)$ in a neighborhood of $(a,v)=(0,0)$. This
solution depends smoothly on $a, \epsilon$ and satisfies
$||v^*(a,\epsilon)||_{H^2}=O(a^2+\epsilon^2)$ if $|a|, \epsilon \geq
0$ small. Using this in $(\ref{eq11})_{(ii)}$ and taking the $L^2$
inner product with $\phi_1(\epsilon)$ yields
\begin{equation}\label{eqBae}
B(a,\epsilon):= -a\lambda_1(\epsilon)||\phi_1(\epsilon)||^2_{L^2}-
\left(N_{ap}(a\phi_1(\epsilon)+v^*)+G(\epsilon),\phi_1(\epsilon)\right)=0
\end{equation}
i.e.
\[
B(a,\epsilon)=-a\lambda_1(\epsilon)||\phi_1(\epsilon)||^2_{L^2}+
\frac{1}{2}\left(f''(u_{ap})\phi_1^2,\phi_1
\right)a^2-(G(\epsilon),\phi_1)+O(a^\mu \epsilon^\nu)=0
\]
as $a,\epsilon \to 0$ with $\mu +\nu \geq 3$. If
$\lambda_1(\epsilon)=d\epsilon +O(\epsilon^2)$ with $d \neq 0$
(indep. of $\epsilon$), then we could apply the implicit function
theorem to $\epsilon^{-2}B(\epsilon \tilde{a},\epsilon)=0$ and find
an $a_*=O(\epsilon)$ satisfying (\ref{eqBae}). However, since by
$(\ref{Lepsilon})_{(i)}$ we have $d=0$, this analysis breaks down.

\appendix
\section{Proof of Lemma \ref{lem1}}\label{apen1}
$u \in L^\infty (\mathbb{R})$ and (H2) imply that $\Delta u \in
L^\infty (\mathbb{R})$ and
\[
\Delta u
(z)=\sum_{k=-\infty}^{\infty}a_k[u(z-k)-u(z)]=\sum_{k=-\infty}^{\infty}a_k\int_{0}^{-k}u'(z+t)dt,\
\ \ z\in \mathbb{R}.
\]Then, by the Cauchy-Schwarz inequality, we get
\[
\left(\Delta u (z)\right)^2\leq
\left(\sum_{k=-\infty}^{\infty}|a_k|\right)\sum_{k=-\infty}^{\infty}|a_k|\left(\int_{0}^{-k}u'(z+t)dt\right)^2\stackrel{(H2)}{\leq}
C\sum_{k=-\infty}^{\infty}|a_k|(-k)\int_{0}^{-k}u'^2(z+t)dt.
\]So,
\[
\int_{-\infty}^{\infty}\left(\Delta u (z)\right)^2dz\leq C
\int_{-\infty}^{\infty}\left(\sum_{k=-\infty}^{\infty}|a_k|(-k)\int_{0}^{-k}u'^2(z+t)dt
\right)dz=
\]
\[
=C\sum_{k=-\infty}^{\infty}|a_k|(-k)\int_{-\infty}^{\infty}\int_{0}^{-k}u'^2(z+t)dtdz
=C\sum_{k=-\infty}^{\infty}|a_k|(-k)\int_{0}^{-k}\int_{-\infty}^{\infty}u'^2(z+t)dzdt=
\]
\[
=C\sum_{k=-\infty}^{\infty}|a_k|k^2||u'||^2_{L^2}\stackrel{(H2)}{\leq}
C ||u'||^2_{L^2}<\infty.
\]
Thus, $\Delta u \in L^2(\mathbb{R}).$

We have
\[
(\Delta
u,u')=\int_{-\infty}^{\infty}u'(z)\sum_{k=-\infty}^{\infty}a_k[u(z-k)-u(z)]dz=
\]
\[
=\sum_{k=-\infty}^{\infty}a_k\int_{-\infty}^{\infty}u'(z)[u(z-k)-u(z)]dz=\sum_{k=-\infty}^{\infty}a_k\int_{-\infty}^{\infty}u'(z+k)[u(z)-u(z+k)]dz=
\]
\[
\stackrel{a_k=a_{-k}}{=}\sum_{k=-\infty}^{\infty}a_k\int_{-\infty}^{\infty}u'(z-k)[u(z)-u(z-k)]dz.
\]
Thus,
\[
2(\Delta
u,u')=\sum_{k=-\infty}^{\infty}a_k\int_{-\infty}^{\infty}\frac{d}{dz}\left\{u(z)u(z-k)-\frac{u^2(z)}{2}-\frac{u^2(z-k)}{2}
\right\}dz\stackrel{(\ref{eqlimits})}{=}0.
\]
\section{Proof of Proposition \ref{prop2}}\label{apen2}
Since zero is a simple eigenvalue of $-L^0$, it follows from regular
perturbation theory (cf. Sec. 14.3 of \cite{chowhale}) that it
perturbs smoothly to a simple eigenvalue $\lambda(\epsilon)$ of
$-L^\epsilon$. The corresponding eigenfunction $\phi (\epsilon)$
with $||\phi(\epsilon)||_{H^2}=1$ also depends smoothly (in the
$H^2$ norm) on $\epsilon \geq 0$ small and
$\phi(0)=\frac{u_0'}{||u_0'||_{H^2}}$. It is easy to show that
$\lambda(\epsilon)$ is the principal eigenvalue of $-L^\epsilon$; we
denote it by $\lambda_1(\epsilon)$ and the corresponding $H^2$
normalized eigenfunction by $\phi_1 (\epsilon)$. (Recall that
$(-L^0\phi,\phi)\geq c||\phi||_{L^2}^2,\ \forall \phi \in
H^2(\mathbb{R}), \ \phi\bot u_0',$ and $\left|\left|\phi_1
(\epsilon)-\frac{u_0'}{||u_0'||_{H^2}} \right|\right|_{H^2}\leq
C\epsilon$, to obtain $(-L^\epsilon \phi,\phi)\geq
c||\phi||_{L^2}^2,\ \forall \phi \in H^2(\mathbb{R}), \ \phi\bot
\phi_1 (\epsilon)$.) We have
\[
\phi_1''+f'(u_{ap})\phi_1-\epsilon \Delta \phi_1=-\lambda_1 \phi_1
\]
and
\[
-\lambda_1 (\phi_1,u_0')
=\left(\phi_1,u_0'''+f'(u_{ap})u_0'-\epsilon \Delta u_0'
\right)=\left(\phi_1,[f'(u_0+\epsilon u_1)-f'(u_0)]u_0'-\epsilon
\Delta u_0' \right).
\]
 Since $\phi_1 (\epsilon)
\stackrel{H^2}{\to}\frac{u_0'}{||u_0'||_{H^2}}$ as $\epsilon \to 0,$
we get
\[-\lim_{\epsilon \to
0}\frac{\lambda_1(\epsilon)}{\epsilon}=\frac{1}{||u_0'||_{L^2}^2}\left(u_0',f''(u_0)u_1u_0'-
\Delta u_0' \right).\] Differentiating (\ref{equ1}) yields
\[
u_1'''+f''(u_0)u_0'u_1+f'(u_0)u_1'=\Delta u_0'
\]
i.e.
\[
\left(u_0',f''(u_0)u_1u_0'- \Delta u_0'
\right)=-\left(u_0',u_1'''+f'(u_0)u_1'
\right)=-\left(u_0'''+f'(u_0)u_0',u_1' \right)=0.
\]
 This and the smoothness of
$\lambda_1(\epsilon)$ gives us $(\ref{Lepsilon})_{(i)}$. From the
smoothness of $\phi_1 (\epsilon)$ (in the $H^2$ norm) we have
$(\ref{Lepsilon})_{(ii)}$.
\\

\emph{Acknowledgments.} The author would like to thank the anonymous referees for carefully reading the paper and for offering several pertinent remarks. 
This work has received funding from the Hellenic Foundation for Research and Innovation (HFRI) and the General Secretariat for Research and Technology (GSRT), under grant agreement No 1889.


\begin{thebibliography}{50}
\selectlanguage{english}

\bibitem[ACR]{aigner}
A.A. Aigner, A.R. Champneys, V.M. Rothos,
\emph{A new barrier to the existence of moving kinks in Frenkel-Kontorova lattices}
Physica D  186,   148--170, (2003).


\bibitem[AFFS]{AFFS} N.D. Alikakos, P.C. Fife, G. Fusco, C. Sourdis, \emph{Analysis of the heteroclinic connection in
a singularly perturbed system arising from the study of crystalline
grain boundaries}, Interfaces Free Bound. 8, no. 2, 159--183,
(2006).
\bibitem[Ar]{arnold} V. Arnold, \emph{Ordinary Differential
Equations}, MIT Press, Cambridge, MA, (1973).
\bibitem[BCC]{bateschenchmaj} P. W. Bates, X. Chen, A. Chmaj,
\emph{Traveling waves of bistable dynamics on a lattice}, SIAM J.
Math. Anal. 35, no. 2, 520--546, (2003).

\bibitem[B]{BZ}P.W. Bates, \emph{On some nonlocal evolution equations arising in materials science}, in Nonlinear dynamics and evolution equations, Amer. Math. Soc., vol 48 of Fields Inst. Commun,  13--52 (2006).

\bibitem[BeK]{BeK} T.I. Belova, A.E. Kudryavtsev, Phys. Usp. 40, 359 (1997).

 \bibitem[BL]{bl}H. Berestycki, P.-L. Lions, \emph{Nonlinear scalar field equations I: Existence of a ground state}, Arch. Rational Mech. Anal. 82, no. 4,
  313--345, (1983).
  \bibitem[BK]{braunK}
  O.M. Braun, Y.S. Kivshar, \emph{The Frenkel-Kontorova Model}, Concepts, Methods,
and Applications, Springer-Verlag, 2004.

\bibitem[CH]{chowhale} S.-N. Chow, J. Hale, \emph{Methods of Bifurcation Theory}, Springer-Verlag,
(1982).

\bibitem[CMS]{chowMallet} S.-N. Chow, J. Mallet-Paret, W. Shen, \emph{Traveling waves in lattice dynamical
systems}, J. Differential Equations 149, no. 2, 248--291, (1998).

\bibitem[dPK]{dPK} M. del Pino and M. Kowalczyk, \emph{Renormalized energy of interacting Ginzburg-Landau
vortex filaments}, J. Lond. Math. Soc. (2) 77, no. 3, 647--665, (2008).


\bibitem[DKY]{dmitriev}S.V. Dmitriev, P.G. Kevrekidis, N. Yoshikawa, \emph{Discrete Klein-Gordon models with static kinks free of the Peierls-Nabarro potential},  J. Phys. A.: Math. Gen. 38, 7617--7627, (2005).

\bibitem[FR1]{rothos}
M. Feckan, V.M. Rothos,    \emph{Kink-like periodic travelling waves for lattice equations with on-site and
intersite potentials}, Dyn. PDEs 2, 357--370, (2005).

\bibitem[FR2]{rothos2}
M. Feckan, V.M. Rothos,
\emph{Travelling waves in Hamiltonian systems on
2D lattices with nearest neighbour interactions},
Nonlinearity 20,  319--341, (2007).

\bibitem[FR3]{rothos3}
M. Feckan, V.M. Rothos,   \emph{Travelling waves of discrete nonlinear Schr\"{o}dinger
equations with nonlocal interactions}, Applicable Analysis  89, 1387--1411, (2010).


\bibitem[HK]{halekocak}J. Hale, H. Kocak, \emph{Dynamics and
Bifurcations}, Springer-Verlag, (1991).
\bibitem[He]{henry}D. Henry, \emph{Geometric Theory of Semilinear Parabolic
Equations}, Lecture Notes in Mathematics Vol. 840, Springer-Verlag,
(1981).
\bibitem[IS]{imry}
Y. Imry, L. Schulman, \emph{Qualitative theory of the nonlinear behavior of coupled
Josephson junctions}, J. Appl. Phys. 49,  749--758, (1978).

\bibitem[IP]{IossPeli}G. Iooss,  D.E. Pelinovsky, \emph{Normal form for travelling kinks in discrete Klein-Gordon lattices},
Physica D 216,    327--345, (2006).

\bibitem[KKCR]{kouk}
 V. Koukouloyannis,
  P.G. Kevrekidis,
 J. Cuevas,      V. Rothos,
\emph{Multibreathers in Klein-Gordon chains with interactions beyond
nearest neighbors}, Physica D {242}, 16--29, (2013).






\bibitem[KZ]{kreiner}
C-F. Kreiner,   J. Zimmer, \emph{Travelling wave solutions for the discrete sine-Gordon equation with
nonlinear pair interaction}, Nonlinear Analysis 70,  3146--3158, (2009).



\bibitem[KT]{kresetrudi}O. Kresse, L. Truskinovsky, \emph{Mobility of lattice
defects: discrete and continuum approaches}, J. Mech. Phys. Solids
51, no. 7, 1305--1332, (2003).

\bibitem[Le]{levi} M. Levi, \emph{Dynamics of discrete Frenkel-Kontorova models}, in: Analysis, et
Cetera: Research Papers Published in Honor of J\"{u}rgen Moser's 60th Birthday,
Academic Press Inc., 1990,  471--494.

\bibitem[MA]{makay}R. MacKay, S. Aubry, \emph{Proof of existence of breathers for time-reversible or hamiltonian
networks of weakly coupled oscillators}, Nonlinearity 7, 1623--1643,
(1994).

\bibitem[OPB]{oxt}O.F. Oxtoby, D.E. Pelinovsky, I.V. Barashenkov, \emph{Travelling kinks in discrete $\phi^4$ models},
 Nonlinearity 19, 217--235, (2006).


\bibitem[SL]{levia} A. Saadatpour, M. Levi, \emph{Traveling waves in chains of pendula}, Physica D 244,  68--73, (2013).




\bibitem[SZE]{savin} A.V. Savin, Y. Zolotaryuk, J.C. Eilbeck,
 \emph{Moving kinks and nanopterons in the nonlinear Klein-Gordon lattice},
 Physica D 138,  267--281, (2000).



\bibitem[S]{Schmidt}
V. H. Schmidt, \emph{Exact solution in the discrete case for solitons propagating in a chain of
harmonically coupled particles lying in double-minimum potential wells}, Phys. Rev. B 20, 4397--4405,
(1979) .

\bibitem[Sc]{scott} A.C. Scott, \emph{Nonlinear Science}, Oxford University Press, (2003).

\bibitem[SZ]{shatah} J. Shatah,  C. Zeng,
\emph{Orbits homoclinic to centre manifolds of conservative
PDEs},
Nonlinearity 16,  591--614, (2003).

\end{thebibliography}
\end{document}